# Exact equations for smoothed Wigner transforms and homogenization of wave propagation

Agissilaos G. Athanassoulis
*Program in Applied and Computational Mathematics, Princeton University, Princeton, NJ-08544, USA*
(April 13 2007)

The Wigner Transform (WT) has been extensively used in the formulation of phase-space models for a variety of wave propagation problems including high-frequency limits, nonlinear and random waves. It is well known that the WT features counterintuitive 'interference terms', which often make computation impractical. In this connection, we propose the use of the smoothed Wigner Transform (SWT), and derive new, exact equations for it, covering a broad class of wave propagation problems. Equations for spectrograms are included as a special case. The 'taming' of the interference terms by the SWT is illustrated, and an asymptotic model for the Schrödinger equation is constructed and numerically verified.

## 1. Introduction

The Wigner Transform (WT), or Wigner function, is a well known object in quantum mechanics [1,20,21,28,31,32,35,38,40], signal processing [5,6,13,14,21,22,23,33], and high-frequency wave propagation [4,8,9,10,11,12,17,18,19,25,26,27,29,30,36,37]. It provides one of the most appropriate ways to describe the passage from quantum to classical mechanics in physics, and at the same time it is close to a proper time-frequency energy density in signal processing. However, its applicability is limited by certain complicated, counterintuitive features, collectively described as interference terms [6,14,16,22,23,24].

Smoothed Wigner Transforms [4,5,6,13], Wigner measures [4,8,17,18,19,25, 26,27,29,30,37,42], Wigner spectra [15,21] and other variants have been proposed as alternatives that keep as much as possible from the WT's good properties while suppressing the features that make it impractical. An area that has used heavily objects of this kind is wave propagation, and more specifically the homogenization of wave equations in various physical contexts. Wigner measures (also called microlocal defect measures) have been developed in order to keep track of nonlinear quantities, e.g. energy, in high-frequency wave propagation, in acoustics, electromagnetism, elasticity and quantum mechanics. WTs and Wigner spectra are also used in the study of high-frequency wave problems with random and/or nonlinear features, including semiconductors and nonlinear optics [1,9,10,11,21,28,36]. Recently, computational schemes for linear, nonlinear, deterministic and stochastic WT-based models have been developed [11,26,27,28].

It seems to be widely recognized however that the numerical propagation of the WT (even with asymptotic, e.g. Liouville, equations) is not a practical approach for high-frequency wave propagation. This is so because of the interference terms; the WT is itself a high-frequency wavefunction in a higher dimensional space than the original wavefunction, and working with it would be more expensive than working with the original problem. Indeed, in all the computational works mentioned above, the subject of the numerical evolution is either a Wigner measure or a smooth phase space density (which may be interesting due to nonlinearities [21]), or a Wigner spectrum (i.e. smoothed in the context of a stochastic formulation [15,21,36]), but never a complete WT corresponding to a realistic high-frequency wavefield. In other words, Wigner homogenization usually consists of two steps: first, one simplifies the equations of motion for the WT – to get for example a Liouville equation – then one uses a reduced representation to keep track of the main features of WT.

A particular approach that has been extensively worked out is to keep a Liouville equation, as far as the equations of motion are concerned, and use Wigner measures to describe asymptotically the WT, e.g. as in [19,25,29]. The Wigner measure is an asymptotic object, with no interference terms and no negative values. The prototype Wigner measure corresponds to a WKB wavefunction $u^\varepsilon(x,t) = A(x,t)e^{\frac{2\pi i}{\varepsilon}S(x,t)}$, $\varepsilon \ll 1$, and is $W^0(x,k,t) = |A(x,t)|^2 \delta(k - \nabla_x S(x,t))$. This approach offers a tractable scheme for asymptotic analysis that does not break down when caustics appear. It has recently been used in conjunction with the level set method in the development of a specialized computational technique [26,27]. The Wigner measure approach thus reclaims computability by solving the $\varepsilon \to 0$ limit problem. This also means that information is lost, in particular around caustics [13].

In this paper we propose a computable asymptotic approach based on smoothed WTs, i.e. WTs convolved with a smooth kernel. The derivation of the exact equations for SWTs, and their asymptotic treatment is presented, along with numerical results confirming their validity and the overall computability of this approach. The results are compared to exact and full numerical solutions of the corresponding PDEs.

**2. The Wigner Transform and the smoothed Wigner Transform**

The Wigner Transform is defined as a sesquilinear mapping,

$$W^\varepsilon : f(x), g(x) \mapsto W^\varepsilon[f,g](x,k) = \int_{y \in \mathbb{R}^n} e^{-2\pi i k y} f\left(x + \frac{\varepsilon y}{2}\right) \overline{g}\left(x - \frac{\varepsilon y}{2}\right) dy.$$

When $f = g$ it is called the Wigner distribution (WD) of $f$, and denoted as $W^\varepsilon[f](x,k)$. The WD of a wavefunction $f(x)$ is a good way to realize a joint breakdown of the wavefunction's energy over space $x$ and wavenumber $k$ - with the caution that it takes on negative values as well.

The WT has a number of important properties, which are relevant both in wave propagation and in signal processing. The books [14,16] are two very important sources, [14] emphasizing the signal processing point of view, and [16] the relation between the WT and pseudodifferential operators. Here we will only mention some properties that we directly use. It can be shown that

$$\int_{k \in \mathbb{R}^n} W^\varepsilon[f](x,k) dk = |f(x)|^2, \quad \int_{x \in \mathbb{R}^n} W^\varepsilon[f](x,k) dx = |\hat{f}(k)|^2.$$

Sometimes this is called the exact marginals property; it motivates the interpretation of the WT as a phase space energy density, since it allows the expression of more familiar forms of the energy in terms of the WT. In fact, even when the $L^2$ norm is not the natural energy of the problem, the amplitude of the wavefield is always a relevant quantity. Moreover, it can be shown that, if $L(x, \varepsilon \partial_x)$ is a pseudodifferential operator with Weyl symbol $L(x,k)$, then

$$\int_{x,k \in \mathbb{R}^n} L(x,k) W^\varepsilon[f,g](x,k) dx dk = \int_{y \in \mathbb{R}^n} \overline{g}(y) L(x, \varepsilon \partial_x) f(y) dy.$$

The last equation, sometimes called the trace formula, allows us to express more general energies, as well as other quantities, like energy flux, in terms of the WT. Any sesquilinear/quadratic observable of the wave propagation can be expressed explicitly in terms of the WT.

We mentioned a few times the interference terms. By 'interference terms' we mean wave patterns that appear in regions of phase space with small or no energy [14,24]. ('Energy', unless otherwise specified, means the $L^2$ norm of the wavefunction, and at the same time the net integral of its WT. Accordingly, 'regions of phase space with significant energy' are regions over which the WT has significant integral). The values of the WD on the interference terms are large, often larger than on the regions that hold significant energy. For real-life signals this results in a complicated, counterintuitive, obscured picture. In other words, the WT exhibits as short waves as the original wavefunction, but in twice the space dimensions. The interference terms are the price we pay for the nice analytic properties of the WT. However, unlike the original wavefunction, we can average the waves out and still keep 'most' of the information.

Cohen's class of distributions [6,7,14,24] is defined as the class of all sesquilinear transforms of the form

$$C[f,g](x,k) = \int_{y,u,z \in \mathbb{R}^n} f\left(u+\frac{y}{2}\right)\bar{g}\left(u-\frac{y}{2}\right) e^{-2\pi i[ky+zx-zu]} \phi_C(z,y)\,dy\,du\,dz$$

for any distribution $\phi_C(z,y)$ (restricting $f,g$ to test functions). An equivalent definition is the class of transforms that results from convolving the WT with a distributional kernel $K_C(x,k) = \mathcal{F}^{-1}_{(z,y) \to (x,k)}[\phi_C(z,y)]$,

$$C[f,g](x,k) = \int_{x',k' \in \mathbb{R}^n} K_C(x-x',k-k')\, W[f,g](x',k')\,dx'\,dk'.$$

Among them there are many attractive alternatives to the WT. When the kernel is chosen to be the WD of some function $h$, $K_C(x,k) = W[h](x,k)$, and if $f=g$, we get the spectrogram of $f$ with window $h$. Although spectrograms have nonnegative values, it is sometimes said that they are too spread out and have smoothed-out many important, relatively fine features [14,24]; this point will emerge in some form in the context of wave propagation as well. It turns out that choices which allow some negative values can be more interesting. Smoothing WTs with an appropriate kernel tames the interference terms, but doesn't necessarily kill them completely; there is a balance between smoothing enough and not smoothing too much. In this connection we define the scaled smoothed Wigner Transform (SWT) as the sesquilinear transform

$$\widetilde{W}^{\sigma_x,\sigma_k;\varepsilon}[f,g](x,k) = \frac{2}{\varepsilon \sigma_x \sigma_k} \int_{x',k' \in \mathbb{R}^n} e^{-\frac{2\pi|x-x'|^2}{\varepsilon \sigma_x^2} - \frac{2\pi|k-k'|^2}{\varepsilon \sigma_k^2}} W^\varepsilon[f,g](x',k')\,dx'\,dk'. \quad (1)$$

This is a WT convolved with a tensor-product Gaussian with space-domain variance proportional to $\varepsilon \sigma_x^2$ and wavenumber-domain variance proportional to $\varepsilon \sigma_k^2$. The scaling is selected with the problem (2) and the high-frequency regime in mind; see Figures 1, 2 for intuition on the scale of smoothing. When the SWT exhibits negative values we say we have subcritical smoothing. When $\sigma_x, \sigma_k$ are chosen so that the SWT coincides with a spectrogram we have critical smoothing, and naturally nonnegative values.

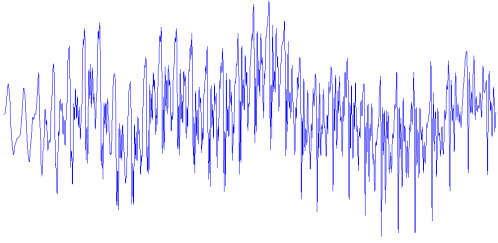

Figure 1a: A segment of human speech

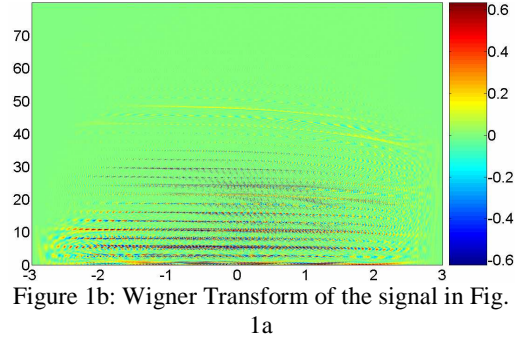

Figure 1b: Wigner Transform of the signal in Fig. 1a

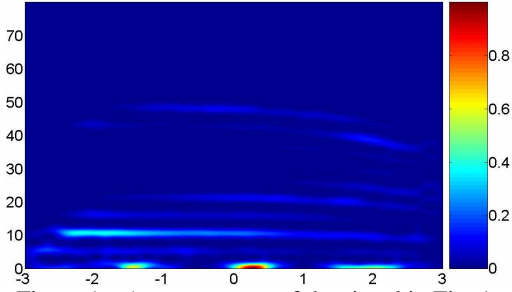

Figure 1c: A spectrogram of the signal in Fig. 1a

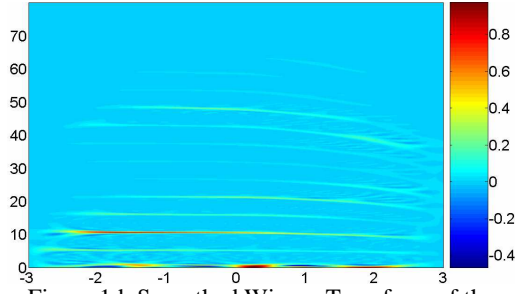

Figure 1d: Smoothed Wigner Transform of the signal in Fig. 1a

## 3. Derivation of the exact equations

Our central result is the derivation of exact equations for the SWT of a wavefunction. The same approach could be used for other transforms from Cohen's class as well [5], but we don't pursue that here. The wave propagation problem we consider is the IVP

$$\varepsilon \frac{\partial}{\partial t} u^\varepsilon (x,t) + L^\varepsilon (x, \varepsilon \partial_x) u^\varepsilon (x,t) = 0, \quad u^\varepsilon (x,0) = u_0^\varepsilon (x) \qquad (2)$$

where $\varepsilon$ is the typical wavelength. The operator's Weyl symbol is assumed to be polynomial (this can be generalized), and may depend on $\varepsilon$, e.g. $L^\varepsilon (x,k) = P(x,k) + \varepsilon Q(x,k)$. (Weyl symbols of differential operators with varying coefficients are typically like that [16,19]). The initial condition $u_0^\varepsilon (x)$ is assumed to be a function with typical wavelength of $O(\varepsilon)$. For now we will also assume it is a test function, but this can be relaxed. The vector (system) counterpart of (2), as well as higher order in time equations, can also be treated; it turns out that the use of the Jordan decomposition is fundamental for those problems. These generalizations, as well as the proofs for non-polynomial symbols, will be published elsewhere. The use of pseudodifferential operators, and specifically of the Weyl calculus, is essential and intimately related with the WT itself [20,35].

The core calculations are the following:

$$\widetilde{W}^{\sigma_x,\sigma_k;\varepsilon}[xf,g](x,k) = \left( x - \frac{\varepsilon}{4\pi i} \frac{\partial}{\partial k} - \frac{\varepsilon \sigma_x^2}{4\pi} \frac{\partial}{\partial x} \right) \widetilde{W}^{\sigma_x,\sigma_k;\varepsilon}[f,g](x,k), \qquad (3a)$$

$$\widetilde{W}^{\sigma_x,\sigma_k;\varepsilon}[f,xg](x,k) = \left( x + \frac{\varepsilon}{4\pi i} \frac{\partial}{\partial k} - \frac{\varepsilon \sigma_x^2}{4\pi} \frac{\partial}{\partial x} \right) \widetilde{W}^{\sigma_x,\sigma_k;\varepsilon}[f,g](x,k), \qquad (3b)$$

$$\widetilde{W}^{\sigma_x,\sigma_k;\varepsilon}\left[\varepsilon\frac{\partial}{\partial x}f,g\right](x,k)=\left(2\pi ik+\frac{\varepsilon}{2}\frac{\partial}{\partial x}+\frac{\varepsilon\sigma_k^2}{2i}\frac{\partial}{\partial k}\right)\widetilde{W}^{\sigma_x,\sigma_k;\varepsilon}[f,g](x,k), \quad (3c)$$

$$\widetilde{W}^{\sigma_x,\sigma_k;\varepsilon}\left[f,\varepsilon\frac{\partial}{\partial x}g\right](x,k)=\left(-2\pi ik+\frac{\varepsilon}{2}\frac{\partial}{\partial x}-\frac{\varepsilon\sigma_k^2}{2i}\frac{\partial}{\partial k}\right)\widetilde{W}^{\sigma_x,\sigma_k;\varepsilon}[f,g](x,k). \quad (3d)$$

With direct use of the identities (3) we can pull out of the SWT any operator with polynomial symbol. A way to write this result compactly is that, for any polynomial $L(x,k)$,

$$\widetilde{W}^{\sigma_x,\sigma_k;\varepsilon}\left[L(x,\varepsilon\partial_x)f,g\right](x,k)=$$
$$=L\left(x-\frac{\varepsilon}{4\pi i}\partial_k-\frac{\varepsilon\sigma_x^2}{4\pi}\partial_x,k+\frac{\varepsilon}{4\pi i}\partial_x-\frac{\varepsilon\sigma_k^2}{4\pi}\partial_k\right)\widetilde{W}^{\sigma_x,\sigma_k;\varepsilon}[f,g](x,k), \quad (4a)$$

$$\widetilde{W}^{\sigma_x,\sigma_k;\varepsilon}\left[f,L(x,\varepsilon\partial_x)g\right](x,k)=$$
$$=L\left(x+\frac{\varepsilon}{4\pi i}\partial_k-\frac{\varepsilon\sigma_x^2}{4\pi}\partial_x,k-\frac{\varepsilon}{4\pi i}\partial_x-\frac{\varepsilon\sigma_k^2}{4\pi}\partial_k\right)\widetilde{W}^{\sigma_x,\sigma_k;\varepsilon}[f,g](x,k). \quad (4b)$$

We discuss the proofs in the Appendix.

We will now use identities (3), (4) to treat the wave propagation problem (2). First of all observe that, if $u^\varepsilon(x,t)$ satisfies (2), due to the sesquilinearity of the SWT it follows that

$$\varepsilon\frac{\partial}{\partial t}\widetilde{W}^{\sigma_x,\sigma_k;\varepsilon}\left[u^\varepsilon(x,t)\right](x,k,t)=$$
$$=\varepsilon\widetilde{W}^{\sigma_x,\sigma_k;\varepsilon}\left[\frac{\partial}{\partial t}u^\varepsilon,u^\varepsilon\right](x,k,t)+\varepsilon\widetilde{W}^{\sigma_x,\sigma_k;\varepsilon}\left[u^\varepsilon,\frac{\partial}{\partial t}u^\varepsilon\right](x,k,t)$$
$$=2\operatorname{Re}\left(\widetilde{W}^{\sigma_x,\sigma_k;\varepsilon}\left[\varepsilon\frac{\partial}{\partial t}u^\varepsilon,u^\varepsilon\right](x,k,t)\right)= \quad (5)$$
$$=2\operatorname{Re}\left(\widetilde{W}^{\sigma_x,\sigma_k;\varepsilon}\left[-L^\varepsilon(x,\varepsilon\partial_x)u^\varepsilon,u^\varepsilon\right](x,k,t)\right).$$

This is the reformulation of (2) that the identities (4) are applicable to. With direct use of equations (4) it follows that (5) can be rewritten as

$$\varepsilon\frac{\partial}{\partial t}\widetilde{W}(x,k,t)+$$
$$+2\operatorname{Re}\left(L^\varepsilon\left(x-\frac{\varepsilon}{4\pi i}\partial_k-\frac{\varepsilon\sigma_x^2}{4\pi}\partial_x,k+\frac{\varepsilon}{4\pi i}\partial_x-\frac{\varepsilon\sigma_k^2}{4\pi}\partial_k\right)\widetilde{W}(x,k,t)\right)=0; \quad (6)$$

where we have used the abbreviation $\widetilde{W}(x,k,t)=\widetilde{W}^{\sigma_x,\sigma_k;\varepsilon}\left[u^\varepsilon(x,t)\right](x,k,t)$. The same equation can be cast in a series. Define

$$A=\frac{1}{4\pi i}\left[\left(\partial_x-i\sigma_k^2\partial_k\right)\partial_y-\left(\partial_k+i\sigma_x^2\partial_x\right)\partial_z\right], \quad (7a)$$

$$B=iA=\frac{1}{4\pi}\left[\left(\partial_x-i\sigma_k^2\partial_k\right)\partial_y-\left(\partial_k+i\sigma_x^2\partial_x\right)\partial_z\right]. \quad (7b)$$

We can take real and imaginary parts of these operators, e.g.

$$\operatorname{Re}(A)=-\frac{1}{4\pi}\left(\sigma_k^2\partial_k\partial_y+\sigma_x^2\partial_x\partial_z\right).$$

Now, for any polynomial symbol $L(x,k)$, equation (6) can be rewritten as

$$\varepsilon \frac{\partial}{\partial t}\widetilde{W}(x,k,t) + \sum_{m=0}^{\infty}\frac{2\varepsilon^m}{m!}\operatorname{Re}[A^m]\operatorname{Re}\left[L(z,y)\big|_{(z,y)=(x,k)}\widetilde{W}(x,k,t)\right] + \\ + \sum_{m=0}^{\infty}\frac{2\varepsilon^m}{m!}\operatorname{Re}[B^m]\operatorname{Im}\left[L(z,y)\big|_{(z,y)=(x,k)}\widetilde{W}(x,k,t)\right]. \quad (8)$$

Equation (8) is the form that can be readily used for formal asymptotic considerations. For its derivation, we Taylor-expand the symbol of equation (6) around $L(x,k)$. Observe also that as long as $L$ is a polynomial, only a finite number of the terms of the series will be nonzero.

**4. High-frequency asymptotics for the Schrödinger equation**

Physically, the high-frequency – or geometrical optics, or semiclassical – asymptotic regime corresponds to wave propagation problems where the coefficients (soundspeed, potential etc) vary on a length scale much larger than the wavelengths that appear. Semiclassical limits of quantum mechanics [18,29,31,32,40], and fluid mechanics [39], are two rich sources of problems in this regime.

Let us see what (5) leads to for the Schrödinger equation with polynomial potential $V(x)$ in the high-frequency regime. Let $u^\varepsilon(x,t)$ satisfy

$$\varepsilon \frac{\partial}{\partial t}u^\varepsilon(x,t) - i\frac{\varepsilon^2}{2}\Delta u^\varepsilon(x,t) + iV(x)u^\varepsilon(x,t) = 0, \quad u^\varepsilon(x,0) = u_0^\varepsilon(x). \quad (9)$$

Then equation (8) for $L(x,k) = -\frac{i}{2}(2\pi k)^2 + iV(x)$ implies

$$\frac{\partial}{\partial t}\widetilde{W}(x,k,t) + \left(2\pi k \frac{\partial}{\partial x} - \frac{V'(x)}{2\pi}\frac{\partial}{\partial k}\right)\widetilde{W}(x,k,t) + \\ + \varepsilon\left(V''(x)\frac{\sigma_x^2}{8\pi^2} + \frac{\sigma_k^2}{2}\right)\frac{\partial^2}{\partial x \partial k}\widetilde{W}(x,k,t) = O(\varepsilon^2), \quad (10) \\ \widetilde{W}(x,k,0) = \widetilde{W}^{\sigma_x,\sigma_k;\varepsilon}[u_0^\varepsilon(x)](x,k)$$

In particular, the leading order of equation (10) is a Liouville equation in phase space. Naturally, the interpretation of the series (8) as a multiple scale expansion is only valid as long as $\sigma_x, \sigma_k = O(1)$, which is of course the default choice following definition (1).

It can be easily seen that making the SWT to coincide with a spectrogram corresponds to setting $\sigma_x \sigma_k = 1$ (in which case the window in the spectrogram is a Gaussian with variance proportional to $\varepsilon \sigma_x^2$). Therefore apparently the Liouville equation is valid for sub-critical (i.e. $\sigma_x, \sigma_k < 1$) as well as critical (i.e. $\sigma_x \sigma_k = 1$) smoothing. In fact this point is a little more subtle – we will come to it shortly.

Both the leading order Liouville equation, and the generalized Fokker-Planck equation obtained by keeping the $O(\varepsilon)$ terms in (10), conserve the total integral of $\widetilde{W}(x,k,t)$, i.e. the $L^2$ norm, in agreement with the respective conservation law for the Schrödinger equation. Indeed, using (10) and integration by parts we find

$$\frac{d}{dt}\int_{x,k}\widetilde{W}(x,k,t)\,dxdk = \int_{x,k}\frac{\partial}{\partial t}\widetilde{W}(x,k,t)\,dxdk =$$

$$= \int_{x,k}\widetilde{W}(x,k,t)\left(\frac{\partial}{\partial x}(2\pi k) - \frac{\partial}{\partial k}\left(\frac{V'(x)}{2\pi}\right) + \varepsilon\frac{\partial^2}{\partial x \partial k}\left(\frac{\sigma_k^2}{2} - V''(x)\frac{\sigma_x^2}{8\pi^2}\right)\right)dxdk + O(\varepsilon^2) =$$

$$= O(\varepsilon^2).$$

In the same way it can be shown that the natural "slow-scale energy" for this problem,

$$\widetilde{\mathcal{E}}(t) = \int_{x,k} L(x,k)\widetilde{W}(x,k,t)\,dxdk, \tag{11}$$

is conserved. The term "slow-scale energy" has the following justification: the quantity

$$\mathcal{E}(t) = \int_{x,k} L(x,k)W(x,k,t)\,dxdk$$

is well known to be equal the natural energy of the problem,

$$\mathcal{E}(t) = \int_{x,k} L(x,k)W(x,k,t)\,dxdk = \langle L(x,\varepsilon\partial_x)u^\varepsilon(x,t), u^\varepsilon(x,t)\rangle;$$

this is a corollary of the so-called trace-formula mentioned in the introduction. It is well known that the two quantities $\mathcal{E}$, $\widetilde{\mathcal{E}}$ are close in high-frequency problems with smooth coefficients [29,33]. A 'smoothed trace formula' is needed here, to allow us to estimate more precisely how close $\int_{x,k} L(x,k)\widetilde{W}(x,k,t)\,dxdk$ is to the actual energy – or to express the energy exactly in terms of $\widetilde{W}(x,k,t)$.

The application of the SWT we focus on in this work is homogenization [19,25,36]. The homogenization scheme proposed here for the Schrödinger equation is to take the SWT of the initial condition of (2), and evolve it in time under the leading order part of (5), typically a Liouville equation for high-frequency problems with smooth coefficients.

This can be used as a slow-scale solution of the problem (2). As we just discussed, moments of WTs describe energy and energy flux, see also [25,36]. The same moments of the corresponding SWTs give a coarser, smoothed out version of the same quantities. Of course the SWT itself is a comprehensive slow-scale representation of the wavefield.

Moreover, knowledge of the SWT can be used in conjunction with time-frequency representations, such as Wilson bases, Malvar-Wilson bases or Gabor frames, to estimate the importance of different atoms (i.e. elements of the basis/frame) for each moment in time. This way a full numerical solution of the problem (2) can be preconditioned.

## 5. Numerical results
*The numerical method*

The WT is computed on a Cartesian grid in phase space with the FFT, with complexity $O(N^2 \log N)$. The complexity for the computation of the SWT is $O(L^2 N^2 \log N)$ where $L$ is the number of sampling points needed for the smoothing kernel. When $\sigma_x, \sigma_k = O(1)$, $L$ is of $O(1)$. Adaptive computation of the SWT that doesn't spend much time on regions of phase space with no energy is also possible [34], and might be essential for two- and higher dimensional problems.

The Liouville equation is solved numerically with the use of particles, i.e. the numerical implementation of the method of characteristics. An initial population of particles is created, so that the SWT can be interpolated up to an error tolerance from its values on them. The trajectory of each particle is computed according to Hamilton's ODEs with a Runge-Kutta solver; the value of the density on each particle remains unchanged in time. The solution at each moment in time is constructed by interpolating the density from its values on the particles.

The results of the method described here are compared with exact solutions, as well as direct numerical solutions of the corresponding Schrödinger equation (9). As we mentioned earlier, being able to compute faster with the SWT, than by direct solution of (9) (or (2), in general), is essential to our point of view. So far, the method we present here appears to be significantly faster than numerical solutions of (9) either with finite-differences or wavelet methods. The difference in complexity appears to be more profound as $\varepsilon \to 0$. Systematic numerical analysis and investigation is in progress, but numerical experiments so far indicate that the time needed to take the SWT of the initial condition and let it evolve it under a Liouville equation (the leading order part of (10) ) is no more than $O\left(\varepsilon^{-2}\log\left(\varepsilon^{-1}\right)\right)$, vs. an observed complexity of at least $O\left(\varepsilon^{-3}\right)$ for standard finite difference methods and $O\left(\varepsilon^{-2.5}\right)$ for a wavelet method.

In all examples we have looked into, the numerical solution of the Liouville equation, i.e. propagation of the particles and interpolation from them, is the slowest part of the numerical solution. In fact, there are a lot of redundant computations there. Typically the number of particles needed to represent the SWT are many-times-over more than enough to represent the flow in phase space described by the Liouville equation – equivalently, Hamilton's ODE. As Ying and Candes point out [41], one can solve the ODE for a sufficient grid of initial conditions, and interpolate for the bulk of the particle population, instead of solving from scratch for each particle. It seems that use of this technique – the phase-flow method, as Ying and Candes call it – would significantly accelerate the method presented here. Let us also note that all the steps of the algorithm, i.e. computation of the SWT, particle propagation, and interpolation, are almost optimally parallelizable.

We also experimented with applying the same computational approach to the WT. It is expected that since it takes more points (bigger support in phase space, but also finer grid needed) it would be accordingly slower. In addition to that, it seems there is an essential stability issue. As the particles move around in phase space, interpolation from them gives very poor results – one could say it becomes very noisy. So, unlike the SWT, a particle population that is sufficient to describe the initial WT, is not necessarily sufficient to evolve it. This kind of behavior shows up whenever the smoothing is too weak.

Three choices for the potential are examined here: a free-space problem, a harmonic oscillator, and a uniform force field. In each case, the corresponding Schrödinger equation is solved, and the SWT and the spectrogram of the initial condition are evolved in time with the corresponding Liouville equation. The $dk$-marginal of each phase-space density is compared to the amplitude of the full solution; the validity of the Liouville equation for the SWT ($\sigma_x, \sigma_k = O(1)$) and its non-validity for the spectrogram are confirmed. One initial condition studied is

$$f^{\varepsilon}(x) = A(x)\exp\left(\frac{2\pi i}{\varepsilon}\left(-\frac{x^4}{4} - x^2 + 2x\right)\right),$$

where $A(x)$ is a smooth envelope defined as

$$A(x) = 0.25\left[\tanh(6.87(x+2.42))+1\right]\left[\tanh(6.87(2.42-x))+1\right].$$

Initial conditions that admit exact solutions are also considered.

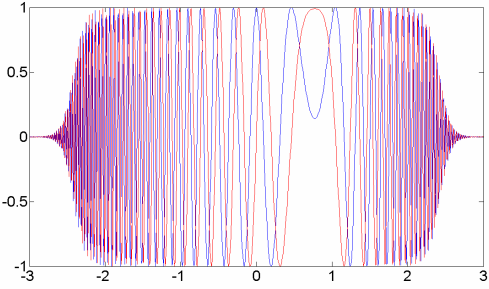

Figure 2a. The function $f^{\varepsilon}(x)$, $\varepsilon = 0.7$ (blue for real part, red for imaginary part).

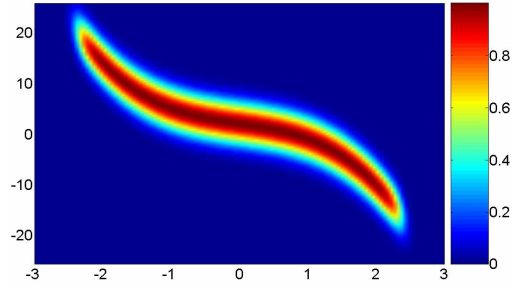

Figure 2b. Spectrogram of $f^{\varepsilon}(x)$, $\varepsilon = 0.7$.

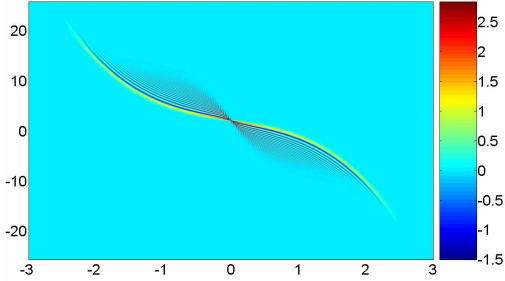

Figure 2c. Wigner Transform of $f^{\varepsilon}(x)$, $\varepsilon = 0.7$.

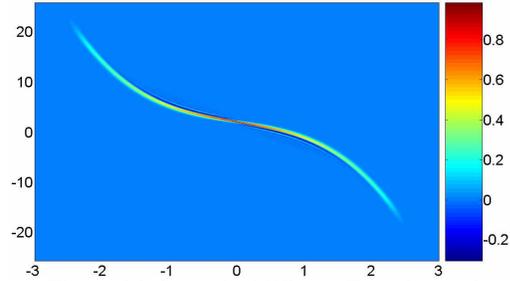

Figure 2d. Smoothed Wigner Transform of $f^{\varepsilon}(x)$, $\varepsilon = 0.7$.

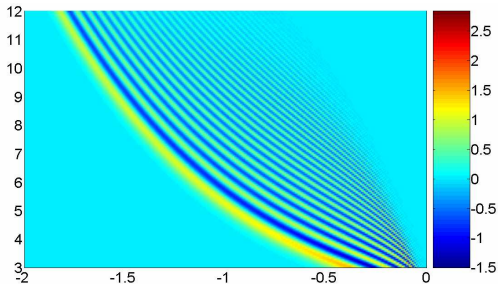

Figure 2e: Interference pattern in the WT (zoomed from Figure 2c)

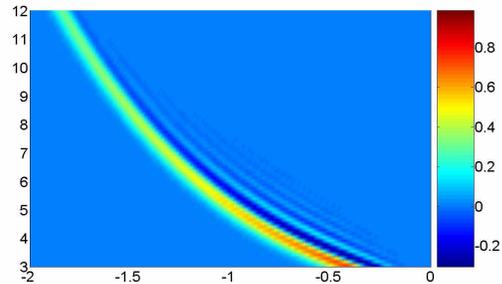

Figure 2f: Interference pattern in the SWT (zoomed from Figure 2d)

*Case studies*
**Case study 1: Free space Schrödinger equation.**
Consider the IVP

$$\varepsilon\frac{\partial}{\partial t}u^{\varepsilon}(x,t) - i\frac{\varepsilon^2}{2}\Delta u^{\varepsilon}(x,t) = 0, \quad u^{\varepsilon}(x,0) = u_0^{\varepsilon}(x). \tag{9a}$$

The corresponding Liouville equation for this problem is

$$\frac{\partial}{\partial t}\widetilde{W}(x,k,t) + 2\pi k\frac{\partial}{\partial x}\widetilde{W}(x,k,t) = 0. \tag{10a}$$

If $u_0^\varepsilon(x)$ is a Gaussian wave packet the exact solution can be computed, i.e. if
$$u_0^\varepsilon(x) = \exp\left[-\left(\mathrm{K}x^2 + \Lambda x + \mathrm{M}\right)\right], \quad \mathrm{K}, \Lambda, \mathrm{M} \in \mathbb{C}, \quad \mathrm{Re}(\mathrm{K}) > 0,$$
then the solution to (6a) is given by
$$u^\varepsilon(x,t) = \frac{1}{\sqrt{2i\varepsilon t \mathrm{K} + 1}} \exp\left[\frac{\Lambda^2 - 4\mathrm{KM} + \dfrac{2i}{\varepsilon t}\left(\mathrm{K}x^2 + \Lambda x + \mathrm{M}\right)}{4\mathrm{K} - \dfrac{2i}{\varepsilon t}}\right].$$

Below is the comparison between the SWT method and the exact solution for (9a) with $\varepsilon = 1$ and initial condition $u_0^\varepsilon(x) = e^{-\frac{1+7i}{0.1}x^2} + e^{-\frac{0.2+3i}{0.1}x^2} + e^{-\frac{0.9-8i}{0.1}x^2}$.

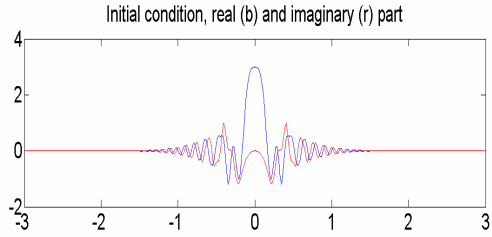

Figure 3a. The initial condition (9)

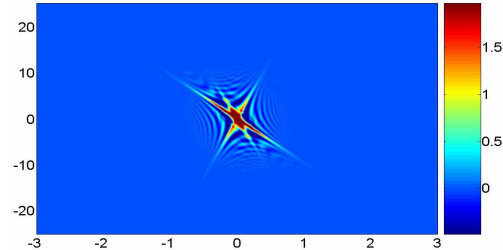

Figure 3b. Wigner Transform of the IC

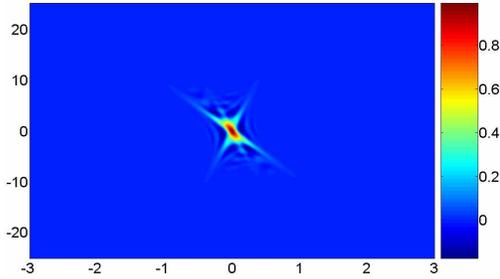

Figure 3c. Smoothed Wigner Transform of the IC

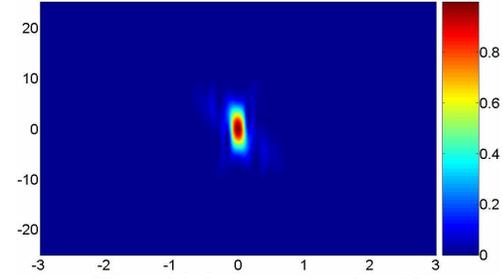

Figure 3d. Spectrogram of the IC

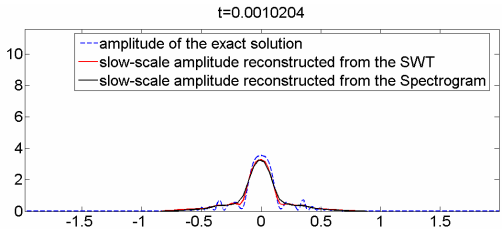

Figure 4a

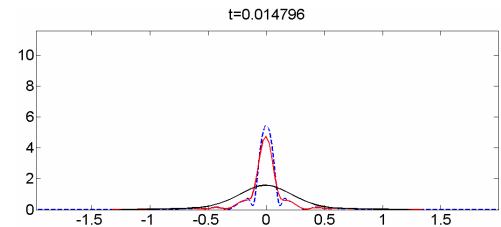

Figure 4b

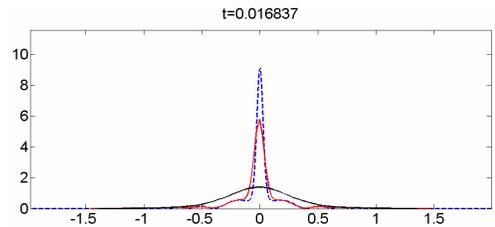

Figure 4c

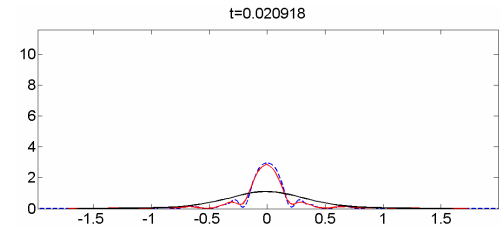

Figure 4d

Figure 4: Snapshots of the evolution in time of the wavefunction's amplitude for problem (6a); exact solution vs SWT and spectrogram evolved in time under eq. (7a); initial condition given in eq. (9).

Consider now equation (9a) with initial condition $u_0^\varepsilon(x) = f^\varepsilon(x)$; $\varepsilon = 0.7$.

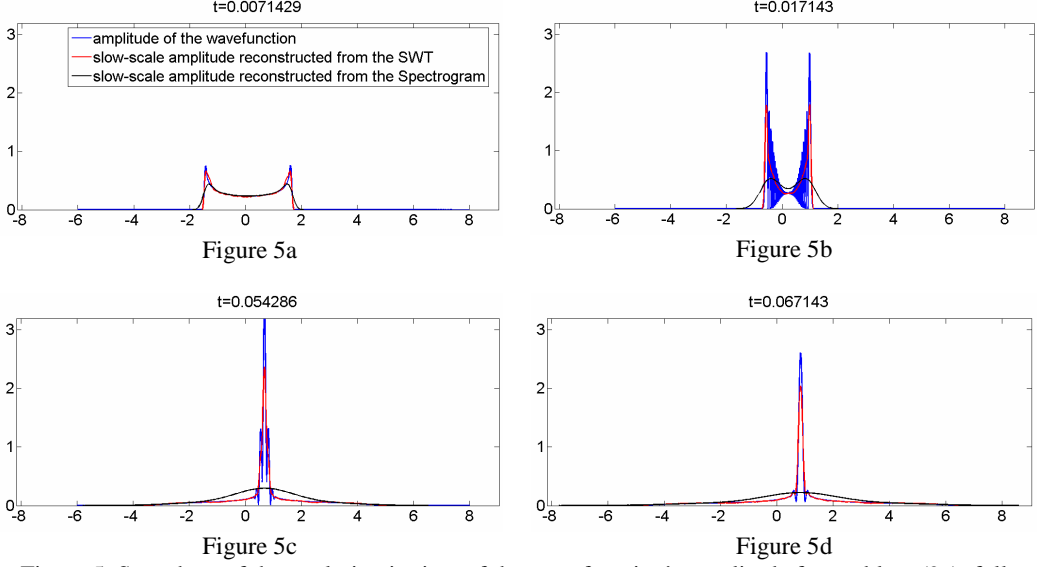

Figure 5a

Figure 5b

Figure 5c

Figure 5d

Figure 5: Snapshots of the evolution in time of the wavefunction's amplitude for problem (9a); full numerical solution vs SWT and spectrogram evolved in time under eq. (10a), initial condition is $u_0^\varepsilon(x) = f^\varepsilon(x)$, $\varepsilon = 0.7$ ( see Fig. 2 for details).

**Case study 2: Quantum harmonic oscillator.**
Consider the IVP
$$\varepsilon \frac{\partial}{\partial t} u^\varepsilon(x,t) - i\frac{\varepsilon^2}{2} \Delta u^\varepsilon(x,t) + i\frac{290}{2} x^2 u^\varepsilon(x,t) = 0, \quad u^\varepsilon(x,0) = u_0^\varepsilon(x). \quad (9b)$$
The corresponding Liouville equation for this problem is
$$\frac{\partial}{\partial t}\widetilde{W}(x,k,t) + 2\pi k \frac{\partial}{\partial x}\widetilde{W}(x,k,t) - \frac{290}{2\pi} x \frac{\partial}{\partial k}\widetilde{W}(x,k,t) = 0. \quad (10b)$$
Time-harmonic solutions for this problem can be constructed in terms of its well-known eigenfunctions. The amplitude of a time-harmonic solution stays constant in time. Again, the SWT leads to satisfactory results and the spectrogram not. Figure 6 corresponds to the ninth eigenfunction.

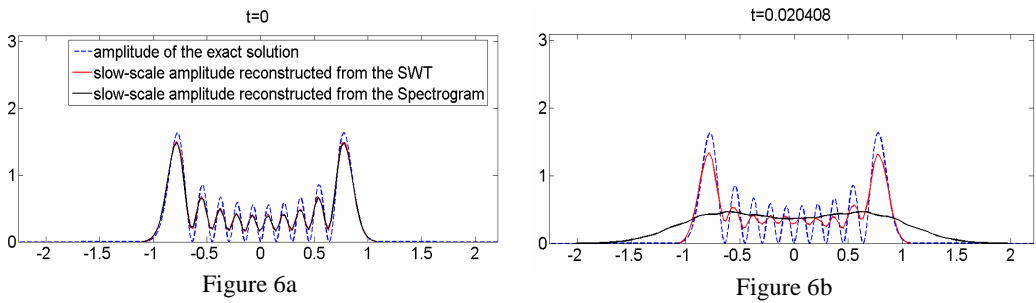

Figure 6a

Figure 6b

Figure 6: Snapshots of the evolution in time of the wavefunction's amplitude for problem (9b); exact solution vs SWT and spectrogram evolved in time under eq. (10b); figure corresponds to the ninth eigenfunction of eq. (9b), $u(x,t) = e^{-i\frac{19}{2}0.7\sqrt{290}\,t} e^{-\frac{1}{2}\frac{\sqrt{290}}{0.7}x^2} H_9\left(x\sqrt{\frac{\sqrt{290}}{0.7}}\right)$, where $H_9(x)$ is the ninth Hermite polynomial. $\varepsilon$ is equal to $0.7$.

Consider now equation (9b) with initial condition $u_0^\varepsilon(x) = f^\varepsilon(x)$, $\varepsilon = 0.7$.

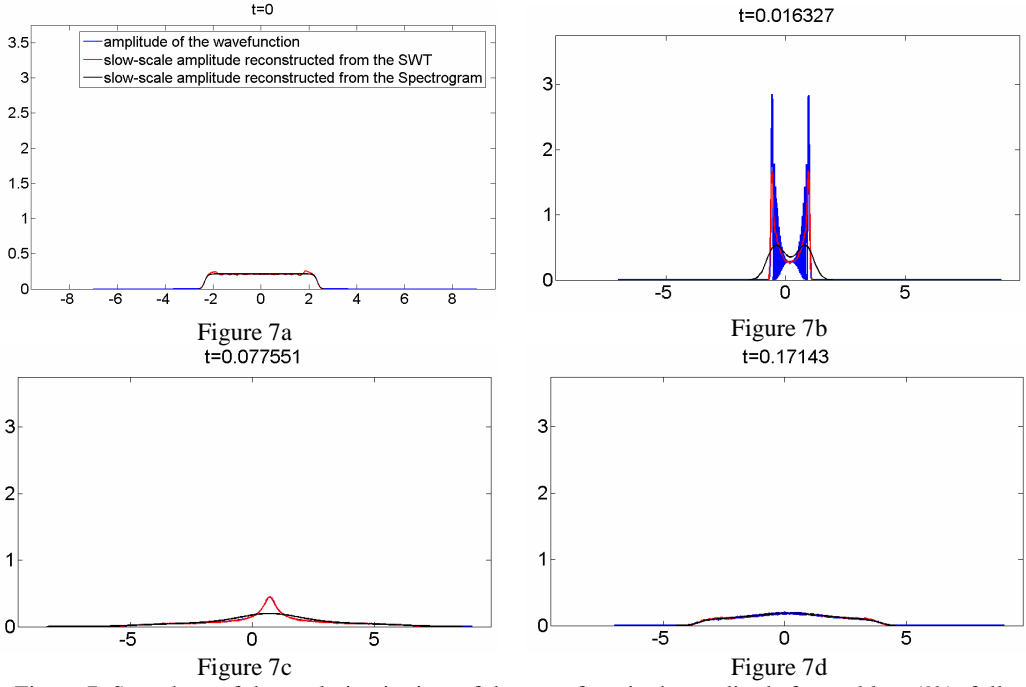

Figure 7: Snapshots of the evolution in time of the wavefunction's amplitude for problem (9b); full numerical solution vs SWT and spectrogram evolved in time under eq. (10b), initial condition is
$u_0^\varepsilon(x) = f^\varepsilon(x)$, $\varepsilon = 0.7$ ( see Fig. 2 for details).

**Case study 3: Uniform force field.**
Consider the IVP

$$\varepsilon \frac{\partial}{\partial t} u^\varepsilon(x,t) - i \frac{\varepsilon^2}{2} \Delta u^\varepsilon(x,t) + i \cdot 2\pi \cdot 300 x \cdot u^\varepsilon(x,t) = 0, \quad u^\varepsilon(x,0) = u_0^\varepsilon(x). \quad (9c)$$

The corresponding Liouville equation for this problem is

$$\frac{\partial}{\partial t}\widetilde{W}(x,k,t) + 2\pi k \frac{\partial}{\partial x}\widetilde{W}(x,k,t) - 300 \frac{\partial}{\partial k}\widetilde{W}(x,k,t) = 0. \quad (10c)$$

Snapshots of the full solution, the SWT based slow-scale amplitude, and the spectrogram based slow-scale amplitude, for $u_0^\varepsilon(x) = f^\varepsilon(x)$, $\varepsilon = 0.7$ follow.

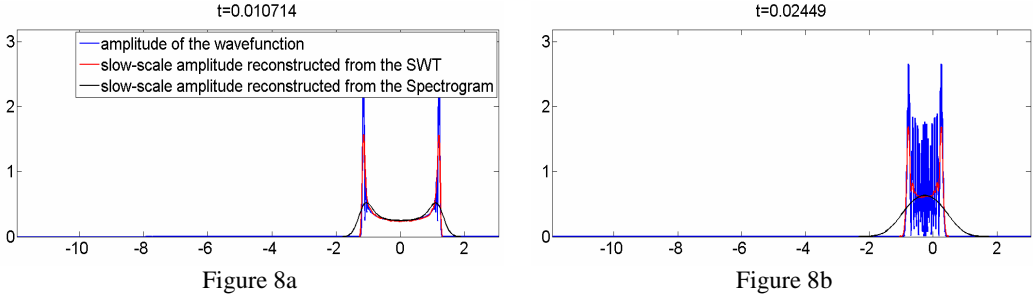

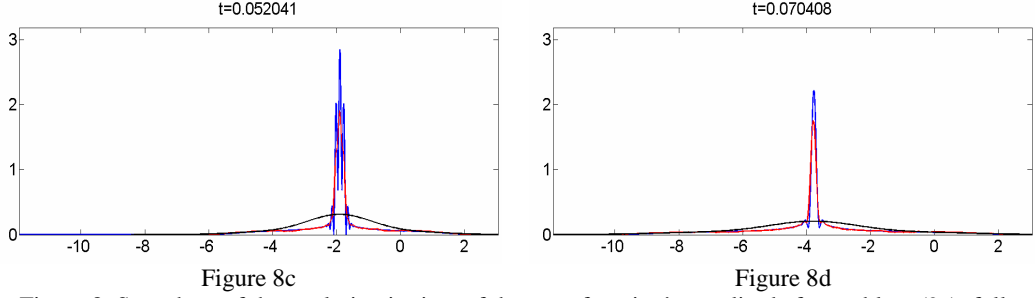

Figure 8c    Figure 8d

Figure 8: Snapshots of the evolution in time of the wavefunction's amplitude for problem (9c); full numerical solution vs SWT and spectrogram evolved in time under eq. (10c), initial condition is $u_0^\varepsilon(x) = f^\varepsilon(x)$, $\varepsilon = 0.7$ ( see Fig. 2 for details).

The full numerical solutions of the Schrödinger equation for the case studies were generously provided by Kostas Politis. They have been carried out using finite-differences and an adaptive, wavelet-based method. The wavelet method is better at handling the development of caustics, as expected. Computations in all cases were done on similar computers, in the MATLAB environment.

### 6. Discussion of the numerical results – Understanding the error

The most striking conclusion of the numerical results is that the method works very well for subcritical smoothing, but fails for critical smoothing (i.e. spectrograms). As we saw earlier, the same formal asymptotics we apply to subcritical smoothing are applicable to the critical case – formally always. Here we will present a very short analysis that helps understand this point; a more complete analysis can be found in [2].

Consider a potential of the form

$$V(x) = ax^s, \text{ where } a \in \mathbb{R} \text{ and } s \in \{0,1,2\}. \tag{12}$$

In any case, it is well known (and can be seen by setting $\sigma_x = \sigma_k = 0$ in (6) or (8) ) the WT of the wavefunction, $W(x,k,t) = W^\varepsilon[u^\varepsilon](x,k,t)$ satisfies (exactly, i.e. not asymtptotically) a Liouville equation. Therefore its evolution can be described in terms of a Hamiltonian flow

$$W(x,k,t) = W(\phi_t(x,k), 0).$$

For these potentials (making use of the fact that Hamiltonian flow happens to be linear) it can be readily seen that, if $\tilde{w}(x,k,t)$ satisfies

$$\frac{\partial}{\partial t}\tilde{w}(x,k,t) + \left(2\pi k \frac{\partial}{\partial x} - \frac{V'(x)}{2\pi}\frac{\partial}{\partial k}\right)\tilde{w}(x,k,t) = 0, \tag{13a}$$

$$\tilde{w}(x,k,0) = \widetilde{W}^{\sigma_x,\sigma_k;\varepsilon}[u_0^\varepsilon](x,k) = \int_{z,y} G_\varepsilon(x-z,k-y) W[u_0^\varepsilon](z,y)\,dzdy, \tag{13b}$$

then

$$\tilde{w}(x,k,t) = \int_{z,y} G_\varepsilon(\phi_t(x-z,k-y)) W(z,y,t)\,dzdy. \tag{14}$$

Put in words, the approximate slow-scale solution that the formal asymptotics lead us to, is not necessarily a very close approximation to $\widetilde{W}^{\sigma_x,\sigma_k;\varepsilon}[u^\varepsilon](x,k,t)$, but rather a slow-scale version of the WT *with smoothing kernel changing with time*. Equation (14) is exact for potentials of the form (12), therefore for the numerical examples presented here. This now allows a better understanding of why spectrograms give "bad" results: when the smoothing kernel is large enough to amount to critical smoothing, its evolution in (14) often dominates the actual propagation, thus giving a very counterintuitive picture. Subcritical smoothing ensures, in many cases, that the evolution of the time-dependent kernel in (14) is small compared to the propagation of the wavefunction.

Equation (14) has been generalized to smooth potentials – there is an additional term accounting for how close the dynamics for the WT are to a Liouville equation (a good measure of how "semiclassical" the problem at hand is). Also, in general the Hamiltonian flows aren't linear, therefore we don't have a convolution, but a more general kernel [2].

**Acknowledgements**

The author would like to thank Profs. G. Papanicolaou, I. Daubechies, G. Makrakis, and L. Ryzhik for helpful discussions. Also, K. Politis for helpful discussions and for providing a part of the numerical results used here. This work has been partially supported by NSF grant DMS-0530865.

**Appendix: Proof of the identities (3)**

For simplicity we present the proof for $\varepsilon = 1$. The scaled version follows in the same lines. First of all note that (3b), (3d) follow from (3a), (3c) respectively due to the sesquilinearity of the SWT. Denote $\phi(z,y) = e^{-\frac{\pi}{2}\left(\sigma_x^2 z^2 + \sigma_k^2 y^2\right)}$. By construction

$$\widetilde{W}[f,g](x,k) = \mathcal{F}_{z,y \to x,k}\left[\phi(z,y) \mathcal{F}^{-1}_{x,k \to z,y}[W[f,g](x,k)]\right], \text{ and}$$

$$W[f,g](x,k) = \mathcal{F}_{z,y \to x,k}\left[\frac{1}{\phi(z,y)} \mathcal{F}^{-1}_{x,k \to z,y}[\widetilde{W}[f,g](x,k)]\right], \text{ or, more compactly,}$$

$$\widetilde{W}[f,g](x,k) = \phi\left(\frac{1}{2\pi i}\partial_x, \frac{1}{2\pi i}\partial_k\right) W[f,g](x,k),$$

$$W[f,g](x,k) = \frac{1}{\phi\left(\frac{1}{2\pi i}\partial_x, \frac{1}{2\pi i}\partial_k\right)} \widetilde{W}[f,g](x,k). \text{ We use the elementary identity}$$

$$W[xf,g](x,k) = \left(x - \frac{1}{4\pi i}\frac{\partial}{\partial k}\right) W[f,g](x,k). \text{ We have}$$

$$\widetilde{W}[xf,g](x,k) = \phi\left(\frac{1}{2\pi i}\partial_x, \frac{1}{2\pi i}\partial_k\right) W[xf,g](x,k) =$$

$$= \phi\left(\frac{1}{2\pi i}\partial_x, \frac{1}{2\pi i}\partial_k\right)\left(x - \frac{1}{4\pi i}\frac{\partial}{\partial k}\right) W[f,g](x,k) =$$

$$= \phi\left(\frac{1}{2\pi i}\partial_x, \frac{1}{2\pi i}\partial_k\right)\left(x - \frac{1}{4\pi i}\frac{\partial}{\partial k}\right) \frac{1}{\phi\left(\frac{1}{2\pi i}\partial_x, \frac{1}{2\pi i}\partial_k\right)} \widetilde{W}[f,g](x,k) =$$

$$= \left(\left[\phi\left(\frac{1}{2\pi i}\partial_x, \frac{1}{2\pi i}\partial_k\right), x - \frac{1}{4\pi i}\frac{\partial}{\partial k}\right] \frac{1}{\phi\left(\frac{1}{2\pi i}\partial_x, \frac{1}{2\pi i}\partial_k\right)} + x - \frac{1}{4\pi i}\frac{\partial}{\partial k}\right) \widetilde{W}[f,g](x,k).$$

The proof of (2a) is completed with the direct computation

$$\left[\phi\left(\frac{1}{2\pi i}\partial_x, \frac{1}{2\pi i}\partial_k\right), x - \frac{1}{4\pi i}\frac{\partial}{\partial k}\right] = \frac{-\sigma_x^2}{4\pi}\frac{\partial}{\partial x}\phi\left(\frac{1}{2\pi i}\partial_x, \frac{1}{2\pi i}\partial_k\right).$$

The same approach, making use of the elementary identity

$$W\left[\frac{\partial}{\partial x}f, g\right](x,k) = \left(2\pi i k + \frac{1}{2}\frac{\partial}{\partial x}\right)W[f,g](x,k),$$

works for (2c). This proof is based on an idea found in [5].